\documentclass{amsproc}
\usepackage{euscript,amstext,amsthm,amssymb}
\usepackage{amsmath}

\theoremstyle{definition}

\theoremstyle{remark}

\numberwithin{equation}{section}

\newcommand{\be}{\begin{equation}}
      \newcommand{\ee}{\end{equation}}
      \newcommand{\ba}{\begin{eqnarray}}
       \newcommand{\ea}{\end{eqnarray}}
\newcommand{\ban}{\begin{eqnarray*}}
       \newcommand{\ean}{\end{eqnarray*}}

\newcommand{\pt}{\partial}

\newcommand{\ra}{\rightarrow}

\newcommand{\bM}{\partial M}

\newtheorem{theo}{Theorem}[section]

\begin{document}
\newtheorem{lem}[theo]{Lemma}
\newtheorem{prop}[theo]{Proposition}
\newtheorem{coro}[theo]{Corollary}

\title{Eta Invariants for Even Dimensional Manifolds}
\author{Xianzhe Dai }
\address{Department of Mathematics, University of  California,
Santa Barbara, California 93106, USA }
\email{dai@math.ucsb.edu}
\thanks{Partially supported by NSF  and NSFC}

\maketitle
\begin{abstract}
In previous work, we introduced eta invariants for even dimensional
manifolds. It plays the same role as the eta invariant of
Atiyah-Patodi-Singer, which is for odd dimensional manifolds. It is
associated to  $K^1$ representatives on  even dimensional manifolds,
and is defined on a finite cylinder, rather than on the manifold
itself. Thus it is an interesting question to find an intrinsic
spectral interpretation of this new invariant. Using adiabatic limit
technique, we give such an intrinsic interpretation.
\end{abstract}

\section{Introduction}


The $\eta$-invariant is introduced by Atiyah-Patodi-Singer in their
seminal series of papers \cite{aps1,aps2,aps3} as the correction
term from the boundary for the index formula on a manifold with
boundary. It is a spectral invariant associated to the natural
geometric operator on the boundary and it vanishes for even
dimensional manifolds (in this case the corresponding manifold with
boundary will have odd dimension). In our previous work with Weiping
Zhang \cite{dz}, we introduced an invariant of eta type for even
dimensional manifolds. It plays the same role as the eta invariant
of Atiyah-Patodi-Singer.

 Any elliptic differential operator on an odd dimensional closed manifold will
have index zero. In this case, the appropriate index to consider is
that of Toeplitz operators. This also fits perfectly with the
interpretation of the index of Dirac operator on even dimensional
manifolds as a pairing between the even $K$-group and $K$-homology.
Thus in the odd dimensional case one considers the odd $K$-group and
odd $K$-homology. For a closed manifold $M$, an element of
$K^{-1}(M)$ can be represented by a differentiable map from $M$ into
the unitary group
$$ g: \ M \longrightarrow { U}(N, {\bf C}), $$
where $N$ is a positive integer.
As we mentioned the appropriate index pairing between the odd $K$-group
and $K$-homology is given by that of the Toeplitz
operator, defined as follows.

Consider $L^2(S(TM)\otimes E)$,  the space of $L^2$ spinor fields\footnote{In this paper,  for
simplicity,  we will generally assume that our manifolds are spin,
although our discussion extends trivially to the case of Clifford
modules} (twisted by an auxilliary vector bundle $E$). It decomposes into an orthogonal direct sum
$$L^2(S(TM)\otimes E)=\bigoplus_{\lambda\in {\rm Spec}(D^E)} E_\lambda,$$
according to the eigenvalues  $\lambda$ of the Dirac operator $D^E$. The ``Hardy space" will be
$$L^2_+(S(TM)\otimes E) =\bigoplus_{\lambda\geq 0}E_\lambda. $$
The corresponding orthogonal projection from $L^2(S(TM)\otimes E)$ to $L^2_+(S(TM)\otimes E)$ will be denoted by $P^E_{\geq 0}$.

The Toeplitz operator $T^E_g$ is then defined as
$$T^E_g=P^E_{\geq 0}gP^E_{\geq 0}:L^2_+\left(S(TM)\otimes E \otimes {\bf C}^N\right)
\longrightarrow L^2_+\left(S(TM)\otimes E \otimes {\bf
C}^N\right).$$ This is a Fredholm operator whose index is given by
\be \label{tifcm} {\rm ind}\, T^E_g=-\left\langle
\widehat{A}(TM){\rm ch}(E){\rm ch}(g),[M]\right\rangle ,\ee where
${\rm ch}(g)$ is the odd Chern character associated to $g$. It is
represented by the differential form (cf. [Z1, Chap. 1])
$$
{\rm ch}(g)= \sum_{n=0}^{\dim M-1\over 2} {n!\over (2n+1)!}{\rm
Tr}\left[\left(g^{-1}dg\right)^{2n+1}\right].
$$

In \cite{dz} we establish an index theorem which generalizes
(\ref{tifcm}) to the case where $M$ is an odd dimensional spin
manifold with boundary $\partial M$. The definition of the Toeplitz
operator now uses Atiyah-Patodi-Singer boundary conditions on
$\partial M$. The self adjoint Atiyah-Patodi-Singer boundary
conditions depend on choices of Lagrangian subspaces $L \subset \ker
D^E_{\partial M}$. We will denote the corresponding boundary
condition by $P_{\partial M}(L)$. The resulting Toeplitz operator
will then be denoted by $T^E_g(L)$.

\begin{theo}[Dai-Zhang] \label{dz} The Toeplitz operator $T^E_g(L)$ is Fredholm with index given
by
$$
 {\rm ind}\, T^E_g(L)=-\left({1\over 2\pi\sqrt{-1}}\right)^{(\dim M+1)/2}\int_M
\widehat{A}\left(R^{TM}\right){\rm
Tr}\left[\exp\left(-R^E\right)\right]{\rm ch}(g) \
$$
$$
-  \overline{\eta}(\partial M, E, g) + \tau_\mu \left(gP_{\partial
M}(L) g^{-1},P_{\partial M}(L) ,  {\mathcal P}_M\right) .
$$
\end{theo}

 Here $\overline{\eta}(\partial M, E, g)$ denotes the
invariant of $\eta$-type for even dimensional manifold $\pt M$ and
the $K^1$ representative $g$. The third term is an interesting new
{\em integer} term here, a triple Maslov index introduced in [KL],
see \cite{dz} for details.

This paper is organized as follows. In Section 2, we review the
definition of the eta invariant for an even dimensional closed
manifold introduced in \cite{dz}. In Section 3, we discuss some
general properties of the invariant. In Section 4, we give an
intrinsic spectral interpretation of the eta invariant. And we end with a conjecture
and a few remarks in the last section.
\newline

\noindent {\bf Acknowledgement}: This is a survey about our previous
joint work, as well as the recent new work with Weiping Zhang. I
would like to thank my collaborator Weiping Zhang for constant
inspiration. Thanks are also due to Matthias Lesch for interesting
conversations.

\section{ An invariant of $\eta$ type for even dimensional
manifolds}

For an even dimensional closed manifold $X$ and a $K^1$
representative $g: X \ra U(N)$, the eta invariant will be defined
in terms of an eta invariant on the cylinder $[0, 1]\times X$ with
appropriate APS boundary conditions.

In general, for a compact manifold $M$ with boundary $\pt M$ with
the product structure near the boundary, the Dirac operator $D^E$
twisted by an hermitian vector bundle $E$ decomposes near the
boundary as \[ D^E = c({\partial \over
\partial x}) ( {\partial \over
\partial x} + D^E_{\partial M} ). \]
The APS projection $P_{\partial M}$ is an elliptic global boundary
condition for $D^E$. However, for self adjoint boundary
conditions, we need to modify it by a Lagrangian subspace of $\ker
D^E_{\partial M}$, namely, a subspace $L$ of $\ker D^E_{\partial M}$
such that $c({\partial \over
\partial x})L=L^\perp \cap (\ker D^E_{\partial M})$. Since $\partial
M$ bounds $M$, by the cobordism invariance of the index, such
Lagrangian subspaces always exist.

The modified APS projection is then obtained by adding the
projection onto the Lagrangian subspace. Let $P_{\partial M}(L)$
denote the orthogonal projection operator from $L^2((S(TM) \otimes
E)|_{\partial M})$ to $L^2_{+}((S(TM)\otimes E)|_{\partial M})\oplus
L$:
$$P_{\partial M}(L) =P_{\partial M}+P_L,$$
where $P_L$ denotes the orthogonal projection  from
$L^2((S(TM)\otimes E)|_{\partial M})$ to $L$.

The pair $(D^E, P^E_{\partial M}(L) )$ forms a self-adjoint elliptic
boundary problem, and $P_{\partial M}(L) $ is called an
Atiyah-Patodi-Singer boundary condition associated to $L$. We will
denote the corresponding elliptic self-adjoint operator by
$D^E_{P_{\partial M}(L) }$.

In \cite{dz}, we originally intend to consider the conjugated
elliptic boundary value problem $D^{E}_{gP_{\partial M}(L)
g^{-1}}$. However, the analysis turns out to be surprisingly
subtle and difficult. To circumvent this difficulty, a
perturbation of the original problem was constructed.

Let $\psi=\psi(x)$ be a cut off function which is identically $1$ in
the $\epsilon$-tubular neighborhood of $\partial M$ ($\epsilon
>0$ sufficiently small) and vanishes outside the $2\epsilon$-tubular
neighborhood of $\partial M$. Consider the Dirac type operator
$$
D^{\psi}=(1-\psi)D^E + \psi gD^E g^{-1}.
$$

The motivation for considering this perturbation is that, near the boundary, the
operator $D^{\psi}$ is actually given by the conjugation of $D^E$,
and therefore, the elliptic boundary problem $(D^{\psi},
gP_{\partial M}(L) g^{-1})$ is now the conjugation of the APS
boundary problem $(D^E, P_{\partial M}(L))$, i.e., this is now
effectively standard APS situation and we have a self adjoint
boundary value problem $(D^{\psi}, gP_{\partial M}(L) g^{-1})$
together with its associated self adjoint elliptic operator
$D^{\psi}_{gP_{\partial M}(L) g^{-1}}$.

The same thing can be said about the conjugation of $D^{\psi}$: \be
\label{pdo} D^{\psi, g}=g^{-1} D^{\psi} g =D^E + (1-\psi)
g^{-1}[D^E, g] . \ee
We will in fact use $D^{\psi, g}$.

We are now ready to construct the eta invariant for even
dimensional manifolds. Given an even dimensional closed spin
manifold $X$, we consider the cylinder $[0, 1]\times X$ with the
product metric. Let $g: \ X \ra U(N)$ be a map from $X$ into the
unitary group which extends trivially to the cylinder. Similarly,
$E \ra X$ is an Hermitian vector bundle which is also extended
trivially to the cylinder. We assume that ${\rm ind}\, D^E_+=0$ on
$X$ which guarantees the existence of the Lagrangian subspaces
$L$.

Consider the analog of $D^{\psi,g}$ as defined in (\ref{pdo}), but
now on the cylinder $[0, 1]\times X$  and denote it by
$D^{\psi,g}_{[0, 1]}$. Here  $\psi=\psi(x)$ is a cut off function on $[0, 1]$ which is identically $1$ for  $0\leq x \leq \epsilon$ ($\epsilon
>0$ sufficiently small) and vanishes when $1-2\epsilon\leq x \leq 1$.  We equip it with the boundary condition
$P_{X}(L) $ on one of the boundary component $\{0\}\times X$ and the
boundary condition ${\rm Id}-g^{-1}P_{X}(L) g $ on the other
boundary component $\{1\}\times X$ (Note that the Lagrangian
subspace $L$ exists by our assumption of vanishing index). Then
$(D^{\psi,g}_{[0, 1]}, P_X(L) , {\rm Id}-g^{-1}P_{X}(L) g)$ forms a
self-adjoint elliptic boundary problem. For simplicity, we will
still denote the corresponding elliptic self-adjoint operator by
$D^{\psi, g}_{[0, 1]}$.

Let $\eta(D^{\psi, g}_{[0, 1]},s)$ be the $\eta$-function of
$D^{\psi, g}_{[0, 1]}$ which, when ${\rm Re}(s)>>0$, is defined by
$$
\eta(D^{\psi, g}_{[0, 1]},s)
=\sum_{\lambda \neq 0}{{\rm sgn}(\lambda)\over
|\lambda|^s},
$$
where $\lambda$ runs through the nonzero
eigenvalues of $D^{\psi, g}_{[0, 1]}$.

By \cite{dw,mu,df}, one knows that the $\eta$-function
$\eta(D^{\psi, g}_{[0, 1]},s)$ admits a meromorphic extension to ${\bf
C}$ with $s=0$ a regular point (and only simple poles). One then
defines, as in \cite{aps1}, the $\eta$-invariant of
$D^{\psi,g}_{[0,1]}$,
$\eta(D^{\psi,g}_{[0, 1]})= \eta(D^{\psi,g}_{[0, 1]},0)$,  and the reduced $\eta$-invariant by
{$$\overline{\eta}(D^{\psi,g}_{[0, 1]})={\dim \ker D^{\psi,g}_{[0, 1]} +
\eta(D^{\psi,g}_{[0, 1]}) \over 2}. $$}

We can also consider the invariant $\overline{\eta}(D^{\psi,g}_{[0,
a]})$, similarly constructed on a cylinder $[0, a] \times X$, and it
turns out to not depend on the radial size of the cylinder $a>0$.
This can be seen by a rescaling argument (cf. [M\"u, Proposition
2.16]).
\newline

\noindent {\bf Definition 2.2.} {We define an invariant of $\eta$
type for the complex vector bundle $E$ on the even dimensional
manifold $X$ (with vanishing index) and the $K^1$ representative
$g$ by
$$
\overline{\eta}(X,E, g)=
 \overline{\eta}(D^{\psi, g}_{[0, 1]}) - {\rm sf} \left\{D^{\psi,g}_{[0, 1]}(s);
 0 \leq s \leq 1 \right\},
$$
where $D^{\psi,g}_{[0, 1]}(s)$ is a path connecting $g^{-1} D^E g$
with $D^{\psi,g}_{[0, 1]}$ defined by
$$
D^{\psi,g}(s)= D^E + (1-s\psi) g^{-1}[D^E, g]  \eqno(2.14)
$$
on $[0, 1]\times X$, with the boundary condition $P_{X}(L) $ on
$\{0\}\times X$ and the boundary condition ${\rm
Id}-g^{-1}P_{X}(L) g$ at $\{1\}\times X$.}
\newline

 It was shown in \cite{dz} that $\overline{\eta}(X, E,
g)$ does not depend on the cut off function $\psi$.

\section{Some properties of the eta invariant}

In this section we look at the properties of our invariant
$\overline{\eta}(X, E, g)$, which depends first of all on the
geometry of the even dimensional manifold $X$, as well as the
hermitian vector bundle $E$ on $X$, and also the $K^1$
representative $g: X \ra U(N)$. An immediate consequence of our
Toeplitz index theorem, Theorem \ref{dz}, is that, when $X=\pt M$ is
the boundary of an odd dimensional compact manifold $M$, the mod
$\mathbb Z$ reduction of $\overline{\eta}(X, E, g)$ is related to
some Chern-Simons invariants:

\be \overline{\eta}(\partial M, E, g) \equiv -\left({1\over
2\pi\sqrt{-1}}\right)^{(\dim M+1)/2}\int_M
\widehat{A}\left(R^{TM}\right){\rm
Tr}\left[\exp\left(-R^E\right)\right]{\rm ch}(g) \ \  {\rm mod} \ \
\mathbb Z. \ee

In this section we will continue to denote the mod $\mathbb Z$
reduction of $\overline{\eta}(X, E, g)$ by the same notation. One
can first study the behavior of $\overline{\eta}(X,E, g)$ under the
metric changes of $X$. Thus let $X_1$ and $X_2$ denote the same
manifold $X$ but with two different Riemannian metrics. Then $X_1 -
X_2= \pt M$ where $M=[0, 1]\times X$ with suitable Riemannian
metrics. Applying the above formula to the current situation yields

\be \overline{\eta}(X_1,E, g) - \overline{\eta}(X_2, E, g) \equiv
-\left({1\over 2\pi\sqrt{-1}}\right)^{(\dim M+1)/2}\int_M
\widehat{A}\left(R^{TM}\right){\rm
Tr}\left[\exp\left(-R^E\right)\right]{\rm ch}(g). \ee

In particular, when $E$ is a flat bundle,

\be \overline{\eta}(X_1,E, g) - \overline{\eta}(X_2, E, g) \equiv -
{\rm rank} E \left({1\over 2\pi\sqrt{-1}}\right)^{(\dim
M+1)/2}\int_M \widehat{A}\left(R^{TM}\right){\rm ch}(g) \ \  {\rm
mod} \ \ \mathbb Z \ee depends only on the rank of the vector bundle
$E$. Thus, if we define

\[ \rho(X, E, g) = \overline{\eta}(X,E, g) - \overline{\eta}(X, {\mathbb C}^{{\rm rank}E}, g),
\ \ {\rm in } \ \mathbb R /\mathbb Z\] then we have deduced that

\begin{theo} The invariant $\rho(X, E, g)$ is independent of the
Riemannian metric on $X$, and hence is an invariant associated to
the manifold $X$, the flat hermitian vector bundle $E$, and the $K^1$
representative $g$. It is a cobordism invariant in the sense that
when $X=\pt M$ is a boundary and $E, \ g$ extends to the interior,
then
\[ \rho(X, E, g) = 0 . \]

\end{theo}

If $E$ is not assumed to be flat but $X=\pt M$ is the boundary of an
odd dimensional compact manifold, the above discussion can be
further refined, as was pointed out in \cite{dz}. We recall it in
the following.

Let $g^{TM}$ resp. $\widetilde{g}^{TM}$, $g^E$ resp.
$\widetilde{g}^E$, and $\nabla^E$ resp. $\widetilde{\nabla}^E$ be
two Riemannian metrics on $M$, two Hermitian metrics on $E$,  and
two connections on $E$. Let $D^E$ resp. $\widetilde{D}^E$ be the
corresponding (twisted) Dirac operators. In order to emphasize the
dependence on the particular geometry of the manifold, we will
denote our eta invariant $\overline{\eta}(\partial M, E, g)$ by the
more explicit notation $\overline{\eta} \left(D^E_{\partial
M},g\right)$

Let $\omega$ be the Chern-Simons form which transgresses the
$\widehat{A}\wedge {\rm ch}$ forms:
$$
d\omega= \left({1\over 2\pi\sqrt{-1}}\right)^{\dim M+1\over 2}
\left(\widehat{A}\left(\widetilde{R}^{TM}\right)\left[\exp\left(-
\widetilde{R}^E\right)\right] -
\widehat{A}\left({R}^{TM}\right)\left[\exp\left(-
{R}^E\right)\right]\right).
$$

Then the following formula describing the variation of
$\overline{\eta} \left(D^E_{\partial M},g\right)$, when
$g^{TM}|_{\partial M}$, $g^E|_{\partial M}$ and $\nabla^E|_{\partial
M}$ change, is proved in \cite{dz}.

$\ $

\begin{theo}[Dai-Zhang] The following identity holds,
$$
\overline{\eta} \left(\widetilde{D}^E_{\partial M},g\right)-
\overline{\eta} \left(D^E_{\partial M},g\right) \equiv-
\int_{\partial M} \omega {\rm ch}(g)\ \  {\rm mod}\ {\mathbb Z}.
$$
\end{theo}

$\ $

One can also study the behavior of the invariant
$\overline{\eta}(X,E, g)$ under the deformations of the $K^1$
representative, as in \cite{dz}. Let $g_t$, $0\leq t\leq 1$ be a
smooth family of $K^1$ representatives $g: X \ra U(N)$. Then
$$\widetilde{\rm ch}\left(g_t,0\leq t\leq 1\right)=
\sum_{n=0}^{(\dim M-1)/2}{n!\over (2n)!}\int_0^1{\rm
Tr}\left[g_t^{-1} {\partial g_t\over\partial
t}\left(g_t^{-1}dg_t\right)^{2n}\right]dt.$$ transgress the odd
Chern character:
$${\rm ch}(g_1)-{\rm ch}(g_0)=d\,\widetilde{\rm ch}\left(g_t,0\leq t\leq 1\right).$$

\begin{theo}[Dai-Zhang] If $\{g_t\}_{0\leq t\leq 1}$ is a smooth family of maps from
$X$ to $U(N)$ and $X$ is a closed even dimensional manifold with
vanishing index, then
$$\overline{\eta} \left(X, E, g_1\right)-\overline{\eta}
\left(X, E,g_0\right)$$
$$\equiv -\left({1\over 2\pi\sqrt{-1}}\right)^{\dim X+1\over 2}\int_{X}
\widehat{A}\left(R^{TX}\right){\rm
Tr}\left[\exp\left(-R^E\right)\right] \widetilde{\rm ch}(g_t,0\leq
t\leq 1)\ \  {\rm mod}\ {\mathbb Z}.$$ In particular, if $g_0={\rm
Id}$, that is, $g=g_1$ is homotopic to the identity map, then \be
\label{wzw} \overline{\eta} \left(X, E,g\right) $$
$$\equiv -\left({1\over 2\pi\sqrt{-1}}\right)^{\dim X+1\over 2}\int_{X}
\widehat{A}\left(R^{TX}\right){\rm
Tr}\left[\exp\left(-R^E\right)\right] \widetilde{\rm ch}(g_t,0\leq
t\leq 1)\ \  {\rm mod}\ {\mathbb Z}.\ee
\end{theo}
\vspace{.3in}

\noindent {\bf Remark} The eta invariant $\overline{\eta}(\bM, g)$
gives an intrinsic interpretation of the Wess-Zumino term in the WZW
theory. When $\bM=S^2$, the Bott periodicity tells us that every
$K^1$ element $g$ on $S^2$ can be deformed to the identity (adding a
trivial bundle if necessary). Hence, (\ref{wzw}) gives another
intrinsic form of the Wess-Zumino term, which is purely local on
$S^2$.

$\ $

Finally, there is also an interesting additivity formula for
$\overline{\eta} (X, E,g)$, as we recall from \cite{dz}.

$\ $

\begin{theo}[Dai-Zhang] Given $f$, $g: X\rightarrow U(N)$, the following identity
holds in ${\mathbb R}/{\mathbb Z}$,
$$
\overline{\eta} \left(X, E,fg\right)= \overline{\eta} \left(X,
E,f\right) + \overline{\eta} \left(X, E,g\right) .
$$
\end{theo}

\section{An intrinsic spectral interpretation}

The usefulness of the eta invariant of Atiyah-Patodi-Singer comes,
at least partially, from the spectral nature of the invariant, i.e.
that it is defined via the spectral data of the Dirac operator on
the (odd dimensional) manifold. Our eta invariant for even
dimensional manifold is defined via the eta invariant on the
corresponding odd dimensional cylinder by imposing APS boundary
conditions. Thus, it will be desirable to have a direct spectral
interpretation in terms of the spectral data of the original
manifold (and the $K^1$ representative). In this section we give
such an interpretation using the adiabatic limit.

First we recall the setup and result from \cite{d}, which is an
extension of \cite{bc} to manifolds with boundary. More precisely,
let \be Y\ra X\stackrel{\pi}{\ra} B \ee be a fibration where the
fiber $Y$ is closed but the base $B$ may  have nonempty boundary.
 Let $g_B$ be a metric on $B$ which is of the
product type near the boundary $\pt B$. Now equip $X$ with a
submersion metric $g$,
\[ g = \pi^{*} g_{B} + g_Y  \]
so that $g$ is also product near $\pt X$. This is equivalent to
requiring $g_Y$ to be independent of the normal variable near $\pt
B$, given by the distance to $\pt B$.

The adiabatic metric $g_x$ on $X$ is given by \be g_x =x^{-2}
\pi^{*} g_{B} + g_Y, \ee where $x$ is a positive parameter.

Associated to these data we have in particular the total Dirac
operator $D^X_x$ on $X$, the boundary Dirac operator $D_x^{\pt X}$
on $\pt X$, and the family of Dirac operators $D_Y$ along the
fibers. If the family $D_Y$ is invertible, then, according to
\cite{bc}, the boundary Dirac operator $D_x^{\pt X}$ is also
invertible for all small $x$,  therefore the eta invariant of $D_x$
with the APS boundary condition,  $\eta(D_x)$, is well-defined.

\begin{theo}  \label{al} Consider the fibration $Y\ra X\ra B$ as above. Assume that the
Dirac family along the fiber, $D_Y$, is  invertible. Consider the
total Dirac operator $D^X_x$ on $X$ with respect to the adiabatic
metric $g_x$ and let $\eta(D^X_x)$ denote the eta invariant of
$D^X_x$ with the APS boundary condition. Then the limit $\lim_{x
\rightarrow 0} \bar{\eta}(D^X_x) = \lim_{x \rightarrow 0}
\frac{1}{2} \eta(D_x)$ exists in ${\Bbb R}$ and
\[\lim_{x \rightarrow 0} \bar{\eta}(D_x) =  \int_B\hat{A}(\frac{R^B}{2\pi})
\wedge \tilde{\eta},  \] where $R^B$ is the curvature of $g_B$,
$\hat{A}$ denote the the $\hat{A}$-polynomial  and $\tilde{\eta}$ is
the $\eta$-form of Bismut-Cheeger \cite{bc}.
\end{theo}

We apply this result to our current situation where $M=[0, 1] \times
X$ fibers over $[0, 1]$ with the fibre $X$. The operator
\[ D^{\psi, g} =D^E + (1-\psi)
g^{-1}[D^E, g]=D^E + (1-\psi) c(g^{-1} dg) \] will be of Dirac type,
and of product type near the boundaries. In order to apply the
adiabatic limit result, we will assume the invertibility condition
that \be \label{ta} \ker [D_X + s\ c(g^{-1} dg)] =0,\ \ \forall\  0
\leq s \leq 1. \ee

Under this assumption there is no spectral flow contribution and
hence \begin{align*} \overline{\eta}(X, g) & = &
 \overline{\eta}(D^{\psi, g}_{[0, 1]}) \\
 & = & \lim_{a\ra \infty} \overline{\eta}(D^{\psi, g}_{[0, a]}) \\
 & = & {\rm the \ adiabatic \ limit \ of \ } \overline{\eta}(D^{\psi, g}_{[0, 1]}) \end{align*}
 is given by the adiabatic limit formula.

By using Theorem \ref{al} we obtain

\begin{theo}[Dai-Zhang] Under the assumption that $\label{ta} \ker [D_X + s\ c(g^{-1} dg)] =0,\ \ \forall\  0
\leq s \leq 1$,
$$
\overline{\eta} \left(X, E,g\right)= \frac{i}{4\pi}
\int_0^1  {\displaystyle \int}_{0}^{\infty} tr_{s}[c(g^{-1} dg)
(D_{X} + s\, c(g^{-1} dg)) e^{-t(D_{X} + s\, c(g^{-1} dg))^2}
 ] dt \ ds .
$$
\end{theo}

For details and further generalization without invertibility
assumption, we refer to \cite{dz2}.

\section{Final remarks}

 Finally we end with a conjecture and some remarks. As we mentioned, the eta type invariant $\overline{\eta}(X, E, g)$, which
we introduced using a cut off function, is in fact independent of
the cut off function. This leads naturally to the question of
whether $\overline{\eta}(X, E, g)$ can actually be defined directly.
The following conjecture is stated in \cite{dz}.

Let $D^{[0, 1]}$ be the Dirac operator on $[0, 1]\times X$. We equip
the boundary condition $gP_{X}(L) g^{-1}$ at $\{0\}\times X$ and the
boundary condition ${\rm Id}-P_{X}(L)$ at $\{ 1 \}\times X$.

Then $(D^{[0, 1]}, gP_{X}(L) g^{-1} , {\rm Id}-P_{X}(L)
 )$ forms a self-adjoint elliptic boundary problem. We denote
the corresponding elliptic self-adjoint operator by
$D^{[0,1]}_{gP_{X}(L) g^{-1} , P_{X}(L) }$.

Let $\eta(D^{[0,1]}_{gP_{X}(L) g^{-1}, P_{X}(L)  },s)$ be the
$\eta$-function of $D^{[0,1]}_{gP_{X}(L) g^{-1} , P_{X}(L)
 }$. By [KL, Theorem 3.1], which goes back to [Gr], one knows
that the $\eta$-function $\eta(D^{[0,1]}_{gP_{X}(L) g^{-1} ,
P_{X}(L)  },s)$ admits a meromorphic extension to ${\bf C}$ with
poles of order at most 2. One then defines, as in [KL, Definition
3.2], the $\eta$-invariant of $D^{[0,1]}_{gP_{X}(L) g^{-1}, P_{X}(L)
}$, denoted by $\eta(D^{[0,1]}_{gP_{X}(L) g^{-1} , P_{X}(L)  })$, to
be the constant term in the Laurent expansion of
$\eta(D^{[0,1]}_{gP_{X}(L) g^{-1}, P_{X}(L)  },s)$ at $s=0$.

Let $\overline{\eta}(D^{[0,1]}_{gP_{X}(L) g^{-1}, P_{X}(L)
 })$ be the associated reduced $\eta$-invariant.

$\ $

\noindent {\bf Conjecture}:
$$
\overline{\eta}(X, E, g)= \overline{\eta}(D^{[0,1]}_{gP_{X}(L)
g^{-1} , P_{X}(L)  }).
$$
\newline

We would also like to say a few words about the technical assumption that
${\rm ind}\, D^E_+=0$ imposed in order to define the eta invariant $\overline{\eta}(X, E, g)$.
The assumption guarantees the existence of the Lagrangian subspaces
$L$ which are used in the boundary conditions. In the Toeplitz index theorem, this
assumption is automatically satisfied since $X=\pt M$ is a boundary. In general, of course,
it may not. However, if one is willing to overlook the integer contribution (as one often does
in applications), this technical issue can be overcome by using another eta invariant, this time
on $S^1 \times X$, as follows. Note that we now have no boundary, hence no need for boundary conditions!

Consider $S^1 \times X= [0, 1] \times X / \sim$ where $\sim$ is the equivalence relation that identifies $0\times X$ with $1\times X$. Let $E_g \ra S^1 \times X$ be the vector bundle which is $E\otimes {\mathbb C}^N$ over $(0, 1) \times X$
and the transition from $0\times X$ to $1\times X$ is given by $g:\ X \ra U(N)$.  Denote by $D_{E_g}$ the Dirac operator  on $S^1 \times X$ twisted by $E_g$.

\begin{prop} One has
\[  \overline{\eta}(X, E, g) \equiv \overline{\eta}(D_{E_g}) \ \ \ {\rm mod} \ \ \mathbb Z . \]
\end{prop}

This is an easy consequence of the so called gluing law for the eta invariant, see \cite{bu,bl,df}.
\newline

\noindent {\bf Remark} It might be interesting to note the duality
that $\overline{\eta} (D^E_{\partial M},g)$
 is a spectral invariant associated to a $K^1$-representative
on an {\it even} dimensional manifold, while the usual
Atiyah-Patodi-Singer $\eta$-invariant ([APS1]) is a spectral
invariant associated to a $K^0$-representative on an {\it odd}
dimensional manifold.

$\ $

Finally, we would like to mention a recent paper of Zizhang Xie
\cite{x} in which he uses our eta invariant to prove, among other
things, an odd index theorem for even dimensional closed manifolds
as well as an odd analog of the relative index pairing formula of
Lesch, Moscovici and Pflaum \cite{lmp}.

\end{document}